\documentclass[12pt]{article}
\usepackage{amsmath, amssymb,amsfonts}
\usepackage{graphicx}
\usepackage{color}
\providecommand{\U}[1]{\protect\rule{.1in}{.1in}}
\usepackage{amssymb}
\usepackage{amsmath, amsfonts}

\newtheorem{teo}{Theorem}

\newtheorem{lemma}{Lemma}

\begin{document}

\title{Inviscid limit for Navier-Stokes equations in domains with permeable boundaries}

\author{N.V. Chemetov* and F. Cipriano** \\
{\small }\\
{\small  * DCM-FFCLRP / University of Sao Paulo, Brazil }\\
{\small E-mail: nvchemetov@gmail.com   }%
\\
{\small and }\\
{\small ** GFM e Dep.\ de Matem{\'a}tica FCT-UNL, Av. Prof. Gama Pinto, 2, 1649-003 Lisboa, Portuga}}
\date{}
\maketitle
\tableofcontents

\begin{abstract}
This work is concerned with 2D-Navier Stokes equations in
a  multiply-connected bounded domain with  permeable walls. The permeability is described by a
Navier type condition. Our aim is to show that the inviscid limit is a
solution of the Euler equations, satisfying the Navier type condition on the
inflow zone of the walls.
\end{abstract}






\section{Introduction and Results}

\label{sec0}

We consider the Navier-Stokes equations for the viscous incompressible
fluids 
\begin{eqnarray}
\mathbf{v}_{t}+\mathrm{div}\,(\mathbf{v}\otimes \mathbf{v})-\bigtriangledown
p &=&{\nu }\Delta \mathbf{v},\quad \mbox{div}\,\mathbf{v}=0\quad \quad \text{%
in}\quad {\Omega _{T}:=\Omega \times (0,T),}  \notag \\
\mathbf{v}(\mathbf{x},0) &=&\mathbf{v}_{0}(\mathbf{x}),\quad \mathbf{x}\in
\Omega ,\quad \text{such that}\quad \mbox{div}\,\mathbf{v}_{0}=0
\label{ICNS}
\end{eqnarray}%
and the Euler equations for the non-viscous incompressible fluids 
\begin{align}
\mathbf{v}_{t}+\mathrm{div}\,(\mathbf{v}\otimes \mathbf{v})-\bigtriangledown
p& =0,\quad \mbox{div}\,\mathbf{v}=0\quad \quad \text{in}\quad {\Omega _{T}},
\notag \\
\mathbf{v}(\mathbf{x},0)& =\mathbf{v}_{0}(\mathbf{x}),\quad \mathbf{x}\in
\Omega \quad \quad \text{with}\quad \mbox{div}\,\mathbf{v}_{0}=0.
\label{ICEE}
\end{align}%
We assume that $\Omega $ is a bounded domain of $\mathbb{R}^{2}$ with the
boundary $\Gamma \in C^{2},$ consisting of the $(n+1)-$connected components $%
S_{0}$, $S_{1}$, \textperiodcentered\ \textperiodcentered\
\textperiodcentered\ , $S_{n}$, and each $S_{k}$ $\ (k=1,...,n)$ lies inside
of $S_{0}$. Here $\mathbf{v}=\mathbf{v}(\mathbf{x},t)$ is the velocity and $%
p=p(\mathbf{x},t)$ is the pressure of the fluid.

Concerning well-posedness, the system \eqref{ICNS} is mostly studied with
the homogeneous Dirichlet boundary condition. Nevertheless, the inviscid
limit with this boundary condition is very difficult to treat and remains an
open problem, essentially due to the creation of a strong boundary layer.
When the system \eqref{ICNS} is supplemented with the so-called Navier
boundary condition, the inviscid limit problem is more treatable.

A correct formulation of the initial-boundary value problem for the Euler
equation \eqref{ICEE} in the planar domain with permeable boundary goes back
to the works of Yudovich \cite{yu}, \cite{yu1}. He proposed the description
of the vorticity on the inlet (for planar flows) and proved the global
existence and uniqueness of the classical solution. In \cite{mor} Morgulis and Yudovich   
investigated stability properties of such flows. In that article  the reader can find 
a clear introduction to the study of the Euler equation with
permeable boundaries.

A remarkable progress has been done in the study of the vanishing viscosity
limit, when the system \eqref{ICNS} is supplemented with Navier boundary
condition for a non-permeable boundary, both in two dimensions (see 
Clopeau et al. \cite{clop}, Kelliher \cite{kel}) 
and in three dimensions (see Beir\~{a}o da Veiga and Crispo \cite{bv-c1}, 
Bellouta et al. \cite{neu}, Xiao and  Xin \cite{x-x}).

The study of the inviscid limit for permeable boundaries is of major
importance for both mathematical and physical communities. Let us refer the
injection/suction systems, created for aircrafts, vehicles, etc, to control
the boundary layers (cf. Marshall \cite{mar}, Schlichting and  Gersten \cite{shu}). In Temam and  Wang \cite{tem} a
conditional result was obtained for Dirichlet's boundary condition. 

Alekseev \cite{alek} studied the vanishing viscosity limit for the steady flows
confined within permeable boundaries, prescribing the vorticity at the
inlet. A non-homogeneous Navier boundary condition was handled in Mucha 
\cite{mucha} for small-sized domains. It was not verified, however, that the
inviscid limit fulfills the boundary condition. 

This article is a version of the article entitled ”Inviscid limit for Navier-Stokes equations in domains with permeable boundaries” published in the journal: Applied Math. Letters'. 33 (2014) 6–11., we refer to \cite{CC5}. After the publication of this paper,  the theory of boundary layers have been developed for Navier-Stokes and Euler equations with Navier slip boundary condition in \cite{C}-\cite{CCC_2}. The Euler equations with sources and sinks has been studied in  \cite{chem}.  

Moreover we have to mention the study of the  non-linear hyperbolic-elliptic systemss, which are similar to the system  of  Euler equations, written  in terms of vorticity - stream function. We refer to the articles \cite{ant2, AC2},  \cite{chem1, chem2}, where the Kruzkov approach  has been generalized  using the kinetic method for the superconductivity model and Keller-Segel's model. Moreover the kinetic method has been developed to porous media models in the articles \cite{chem3}-\cite{chem5}.  The stochastic perturbation of  these hyperbolic-elliptic systems have been studied in the articles \cite{ACC21}, \cite{CC0}.

In the present work we  tackle
the issue of the inviscid limit for permeable boundaries. We establish the
convergence up to the boundary, as the viscosity $\nu \rightarrow 0^{+},$ of
the solutions of the Navier-Stokes equations with the non-homogeneous Navier
type boundary condition%
\begin{equation}
\mathbf{v}\cdot \mathsf{n}=a\quad \mbox{ on }\Gamma _{T},  \label{eq1.39}
\end{equation}%
\begin{equation}
rot(\mathbf{v})=b\quad \mbox{ on }\Gamma _{T},  \label{eqC2}
\end{equation}%
to solutions of the Euler equations satisfying the same boundary condition
on the injection zone of the boundary. Here $\mathsf{n}=(n_{1},n_{2})$
corresponds to the unit outward normal to the boundary $\Gamma .$ The
function $a$ represents the quantity of the inflow and outflow fluid through 
$\Gamma $ with $\int_{\Gamma }a(\mathbf{x},t)\,\,d\mathbf{x}=0,$ $\ t\in
\lbrack 0,T],$ and $\mathrm{rot}(\mathbf{v})=\partial _{x_{1}}v_{2}-\partial
_{x_{2}}v_{1}$ denotes the vorticity of $\mathbf{v}=(v_{1},v_{2})$. The
Navier-Stokes system 
\eqref{ICNS}, \eqref{eq1.39}-\eqref{eqC2} ($NSS$)\ is well posed.

Let us recall that the Navier-Stokes equations \eqref{ICNS} and the Euler
equations \eqref{ICEE} \ can be written in terms of the vorticity $\omega :=%
\mathrm{rot}(\mathbf{v})$ and the velocity $\mathbf{v}$ by 
\begin{eqnarray}
\partial _{t}{}\omega +\mathrm{div}\left( \mathbf{v\omega }\right) &=&{\nu }%
\Delta \omega ,\quad \quad \mbox{div}\mathbf{v}=0\;\quad \quad \mbox{ in }%
\;\Omega _{T},  \label{nsw} \\
\omega (\mathbf{x},0) &=&\omega _{0}(\mathbf{x}),\quad \mathbf{x}\in \Omega
\quad \quad \text{with}\quad \omega _{0}:=\mathrm{rot}\left( \mathbf{v}%
_{0}\right) .  \label{IC-nsw}
\end{eqnarray}%
and%
\begin{eqnarray}
\partial _{t}\omega +\mathbf{v\cdot }\nabla \omega &=&0,\quad \quad %
\mbox{div}\mathbf{v}=0\;\quad \quad \text{in }\Omega _{T},  \label{eq6} \\
\omega (\mathbf{x},0) &=&\omega _{0}(\mathbf{x}),\quad \mathbf{x}\in \Omega ,
\label{7eq7}
\end{eqnarray}%
respectively. The trajectories of the particles for the hyperbolic transport
equation \eqref{eq6} start at the initial moment $t=0$ and on the inflow
region $\Gamma _{T}^{-}.$ We use the notations $\Gamma _{T}^{i}=\cup
_{0<t<T}\Gamma ^{i}(t),$ $i=-,+$ and $0,$\ \ for the parts of $\Gamma _{T},$
where $a<0,$ $a>0$ and $a=0,$ respectively.\ \ Therefore the boundary
condition for the Euler equations \eqref{eq6}-\eqref{7eq7} (or \eqref{ICEE})
should be imposed only on inflow region%
\begin{equation}
\omega =b,\quad \mbox{ on }\Gamma _{T}^{-}.  \label{ome1}
\end{equation}%
We end this short introduction with the statement of our result. The
remaining part of the article is devoted to its proof.

\begin{teo}
\label{teo1} We assume that the data $\,$\thinspace $\mathbf{v}_{0},$\
\thinspace $a,$ $b\,\,\,$\ satisfy the regularity conditions%
\begin{eqnarray}
\mathbf{v}_{0} &\in &W_{p}^{1}(\Omega ),\;\;\; a\in L_{4}(0,T;W_{p}^{1-\frac{1}{p}}(\Gamma ))\cap W_{2}^{1}(0,T;W_{2}^{-\frac{1}{2}%
}(\Gamma )),  \notag \\
b &\in &L_{2}(0,T;W_{p}^{1-\frac{1}{p}}(\Gamma )\cap W_{2}^{1}(0,T;W_{p}^{-%
\frac{1}{p}}(\Gamma ))\;\text{ }\;\text{ with }\;p\in (2,+\infty ).
\label{eq00sec12}
\end{eqnarray}%
There exists a subsequence of $\{\omega _{\nu },\mathbf{v}_{\nu }\},$ being
the solutions of ($NSS$), such that 
\begin{align}
\omega _{{\nu }}& \rightharpoonup \omega \quad \mbox{ weakly}-\ast 
\mbox{ in
}L_{\infty }(0,T;L_{p}(\Omega )),\qquad \nu \nabla \omega _{\nu} \rightharpoonup 0, \;\; \nu \nabla \mathbf{v}_{\nu }
 \rightharpoonup 0\quad \mbox{ weakly in }L_{2}(\Omega _{T}),  \notag \\
\mathbf{v}_{\nu }& \rightharpoonup \mathbf{v}\quad \mbox{ weakly in
}L_{\infty }(0,T,\,W_{p}^{1}(\Omega )) \cap W^1_2 (\Omega _T ) \; \mbox{ and strongly  in
} \; C(0,T,\,C^{\alpha }(\overline{\Omega }))\mbox{ for }\alpha <1-2/p
\label{transition}
\end{align}%
and $\{\omega ,\mathbf{v}\}$ is a solution of the Euler system
\eqref{eq1.39}, \eqref{eq6}--\eqref{ome1} 
 in the following sense%
\begin{align}
\int_{\Omega _{T}}\omega (\psi _{t}+\mathbf{v}\,\cdot \nabla \psi )\;d%
\mathbf{x}dt+\int_{\Omega }\omega _{0}\;\psi (\mathbf{x},0)\;d\mathbf{x}&
=\int_{\Gamma _{T}^{-}}\,a\,b\,\psi \;d\mathbf{x}dt,  \label{eq7} \\
\mathrm{rot}(\mathbf{v})& =\omega ,\;\quad \mathrm{div}(\mathbf{v})=0\quad 
\text{ a.e. in }\,\Omega _{T},  \label{eq8} \\
\mathbf{v}\cdot \mathsf{n}& =a\quad \text{ a.e. on }\,\Gamma _{T},
\label{eq10}
\end{align}%
for every test function $\psi \in C^{1,1}(\overline{\Omega }_{T})$ with 
\begin{equation}
\mbox{ supp}\,(\psi )\subset \left( \Omega \times \left[ 0,T\right) \right)
\cup \Gamma _{T}^{-}.  \label{psi1}
\end{equation}
Moreover there exists $p \in  W^1_2  (\Omega_T )$,
 such that the pair $\{\mathbf{v} ,\,  p \}$ fulfills  the Euler equations
\eqref{ICEE} a.e. in $\Omega_T $. \end{teo}

\section{Proof of Theorem\protect\ref{teo1}}

\label{sec45}

The proof of Theorem \ref{teo1} is a consequence of three lemmas, given in
this section. 



But first let us remember some useful results for the
divergence free functions, defined in multiply-connected domains. For more details about it we refer to Yudovich \cite{yu}, \cite{yu1}. 
Let us consider the solutions $h_{{k}}(\mathbf{x})$\ and 
$h_{\nu }(\mathbf{x},t),$ $h_{a}(\mathbf{x},t)$ of the systems
\begin{equation}
\left\{ 
\begin{array}{l}
-\Delta h_{{k}}=0\quad \quad \text{in }\Omega , \\ 
h_{{k}}=0\quad \text{on }S_{0},\quad h_{{k}}=\delta _{ik}\quad \text{on }%
S_{i},%
\end{array}%
\right. \qquad \left\{ 
\begin{array}{l}
-\Delta h_{{\nu }}=\omega _{{\nu }}\quad \text{ in }\Omega , \\ 
h_{{\nu }}=0\quad \text{ on }\Gamma ,%
\end{array}%
\right. \qquad \left\{ 
\begin{array}{l}
-\Delta h_{{a}}=0\quad \text{ in }\Omega , \\ 
\frac{\partial h_{{a}}}{\partial \mathbf{n}}=a\quad \text{ on }\Gamma%
\end{array}%
\right.  \label{eq13''sec3}
\end{equation}%
with $\ i,k=1,...,n\ \ $and $\ \int_{\Omega }h_{{a}}~d\mathbf{x}=0\ \ \ $%
a.e. on\quad $(0,T).$ \ 

Let us denote $\mathbf{u}_{k}:=\nabla ^{\perp }h_{k}=(\partial
_{x_{2}}h_{k},-\partial _{x_{1}}h_{k})$, $\mathbf{u}_{{\nu }}:=\nabla
^{\perp }h_{\nu },$\ $\mathbf{a}:=\nabla h_{a}$\ and define the matrix $%
\mathbf{A}=(a_{ik})$ with$~a_{ik}:=\int_{\Omega }\mathbf{u}_{{i}}\mathbf{u}_{%
{k}}~d\mathbf{x}.$ It is not difficult to show that $\mathbf{u}_{1},\dots ,%
\mathbf{u}_{n}$ are linearly independent, then the Gram matrix $\mathbf{A}$
is positive definite, that is
\begin{equation}
det(\mathbf{A}) \neq 0\quad \text{and\quad }  \exists c>0: 
\quad  c|\boldsymbol{\lambda }|^{2}\leqslant (\mathbf{A}\boldsymbol{\lambda },%
\boldsymbol{\lambda }),\quad \forall \boldsymbol{\lambda }\in \mathbb{R}%
^{n}.  \label{inv}
\end{equation}%
The solution $\mathbf{v}_{{\nu }}$ of 
($NSS$), satisfying \eqref{eq8}-\eqref{eq10}, can
be expressed in the form%
\begin{equation}
\mathbf{v}_{{\nu }}=\mathbf{u}_{{\nu }}+\sum\limits_{k=1}^{n}\lambda _{k}(t)%
\mathbf{u}_{{k}}+\mathbf{a}\quad \quad \text{in\quad }\Omega _{T}.
\label{VVV}
\end{equation}%
The direct calculations show that $\boldsymbol{\lambda }(t)=(\lambda
_{1}(t),...,\lambda _{n}(t))^{T}$\ is the \textit{unique} solution of the
linear system%
\begin{equation*}
\mathbf{A}\boldsymbol{\lambda }=\mathbf{f\quad }\text{with}\mathbf{\quad f}%
=(f_{1},...,f_{n})^{T},~\quad f_{k}(t):=\int_{S_{k}}(\mathbf{v}_{{\nu }}-%
\mathbf{a})\cdot \mathsf{s}~d\mathbf{x}+\int_{\Omega }\omega _{{\nu }}h_{{k}%
}~d\mathbf{x},
\end{equation*}
where $\mathsf{s}:=(n_{2},-n_{1})$ is the tangent
vector to $\Gamma $. 

The functions $\mathbf{u}_{k}$ $\ (k=1,...,n)$, $\mathbf{u}_{{\nu }}$ and $%
\mathbf{a}$ satisfy the Calderon-Zygmund%
\'{}%
s estimates (Theorems 1.8, 1.10 on pages 12, 15 and Proposition 1.2,
p. 14 in Girault and  Raviart \cite{gir}). \ In addition,
using the compact embedding theorem 
\begin{equation}
W_{p}^{1}(\Omega )\hookrightarrow C^{\alpha }(\overline{\Omega })\qquad 
\text{with\quad }\alpha < 1-2/p,  \label{emb}
\end{equation}%
we have 
\begin{equation*}
||\mathbf{u}_{{k}}||_{C(\overline{\Omega })}\leqslant C||\mathbf{u}_{{k}%
}||_{W_{p}^{1}(\Omega )}\leqslant C_{p},\quad \quad ||\mathbf{u}_{{\nu }%
}||_{W_{p}^{1}(\Omega )}\leqslant C_{p}||\omega _{{\nu }}||_{L_{p}(\Omega )},
\end{equation*}%
\begin{equation}
||\mathbf{a}||_{C(\overline{\Omega })}\leqslant C||\mathbf{a}%
||_{W_{p}^{1}(\Omega )}\leqslant C_{p}||a||_{W_{p}^{1-\frac{1}{p}}(\Gamma )},\quad ||\partial _{t}\mathbf{a}%
||_{L_{2}(\Omega _{T})}\leqslant C||\partial _{t}a||_{L_{2}(0,T;W_{2}^{-%
\frac{1}{2}}(\Gamma ))}
\label{ell-1}
\end{equation}%
with the constants $C_{p}$ depending on $p<\infty .$ 
From \eqref{VVV}, \eqref{ell-1} we obtain for a.a. $t\in (0,T)$%
\begin{equation}
||(\mathbf{v}_{{\nu }}-\mathbf{a})(\cdot ,t)||_{C(\overline{\Omega }%
)}\leqslant C\Vert (\mathbf{v}_{{\nu }}-\mathbf{a})(\cdot ,t)\Vert
_{W_{p}^{1}(\Omega )}\leqslant C(\Vert \omega _{{\nu }}(\cdot ,t)\Vert
_{L_{p}(\Omega )}+|\boldsymbol{\lambda }(t)|).  \label{emb1}
\end{equation}

\subsection{Estimates independent of $\protect\nu $}
\label{sec452}

In this subsection we obtain estimates for $\mathbf{v}_{\nu }$ and $\omega
_{\nu }$ independently of the viscosity $\nu $.  In the first
Lemma we estimate the $L_{2}-$ norm of the velocity $\mathbf{v}_{\nu }$ by
the $L_{p}-$norm of the corresponding vorticity $\omega _{\nu }.$
\begin{lemma}
\label{lem1} Assume that the hypothesis (\ref{eq00sec12}) hold, then there
exists a unique solution $\{\omega _{\nu },\mathbf{v}_{\nu }\}$ of ($NSS$),
satisfying the estimate 
\begin{equation}
\Vert (\mathbf{v}_{\nu }-\mathbf{a})(\cdot ,t)\Vert _{L_{2}(\Omega
)}^{2} + \nu \int_{0}^{t} \Vert D (\mathbf{v}_{\nu }-\mathbf{a})(\cdot ,r)\Vert _{L_{2}(\Omega 
)}^{2} \, dr \leqslant C\left( \int_{0}^{t}f(r)\Vert \omega _{\nu }(\cdot ,r)\Vert
_{L_{p}(\Omega )}^{2}\,dr+1\right)  \label{L2estimate}
\end{equation}%
for any $t\in \lbrack 0,T]$, where $f(t)\in L_{1}(0,T)$ depends only on the
data $\mathbf{v}_{0},$ $a,$ $b$. The constant $C$ is independent of $\nu .$
\end{lemma}

\textbf{Proof}. Since $\Gamma \in C^{2},$ if we parametrize the boundary $\Gamma $ \ by the arc
length $s$, the curvature $k=k(s)$ of $\Gamma $ is a continuous function.
 Applying Lemma 1 of Chemetov and  Antontsev
  \cite{C} the Navier type slip
boundary condition \eqref{eqC2} is equivalent to the following boundary
condition%
\begin{equation}
2D(\mathbf{v})\mathsf{n}\cdot \mathsf{s}+2k\mathbf{v}\cdot {\mathsf{s}}%
=g\quad \mbox{ on }\Gamma _{T},  \label{oq}
\end{equation}%
where $D(\mathbf{v}):=\frac{1}{2}[\nabla \mathbf{v}+(\nabla \mathbf{v})^{T}]$
is the rate-of-strain tensor and $g:=b+2\frac{\partial a}{\partial s}.$ 

It is easy to check that $\mathbf{z}_{\nu }:=\mathbf{v}_{\nu }-\mathbf{a}$
satisfies the system 
\begin{equation}
\begin{cases}
\partial _{t}\mathbf{z}_{\nu }+\mathbf{a}\cdot \nabla \mathbf{z}_{\nu
}-\bigtriangledown p_\nu ={\nu }\Delta \mathbf{z}_{\nu }+F_{\nu },\text{\quad
\quad }\mathrm{div}\,\mathbf{z}_{\nu }=0\text{\quad }\;\text{ in}\;\Omega
_{T}, \\ 
\mathbf{z}_{\nu }\cdot \mathsf{n}=0,\text{\quad \quad }2D(\mathbf{z}_{\nu
})\,\mathsf{n}\cdot \mathsf{s}+2k\mathbf{z}_{\nu }\cdot {\mathbf{\mathsf{s}}}%
=\widetilde{g}\text{\quad \quad on}\;\Gamma _{T}, \\ 
\mathbf{z}_{\nu }(0,\mathbf{x})=\mathbf{v}_{0}(\mathbf{x})-\mathbf{a}(0,%
\mathbf{x})\text{\quad in}\;\Omega%
\end{cases}
\label{NSVV}
\end{equation}%
with $F_{\nu }:={\nu }\Delta \mathbf{a}-\partial _{t}\mathbf{a}-\mathbf{z}%
_{\nu }\cdot \nabla \mathbf{z}_{\nu }-\mathbf{a}\cdot \nabla \mathbf{a-z}%
_{\nu }\cdot \nabla \mathbf{a}$ and $\widetilde{g}:=g-\left[ 2D(\mathbf{a})\,%
\mathsf{n}\cdot \mathbf{\mathsf{s}}+2k\mathbf{a}\cdot {\mathbf{\mathsf{s}}}%
\right] $. The
solvability and the uniqueness of \ \eqref{NSVV} 
are standard results (see Kelliher \cite{kel}).

Multiplying the first equation in \eqref{NSVV} by $\mathbf{z}_{\nu }$,
integrating over $\Omega $ and using \eqref{emb1}, we obtain%
\begin{equation}
\frac{1}{2}\frac{d}{dt}||\mathbf{z}_{\nu }\Vert _{L_{2}(\Omega )}^{2}+2\nu
\int_{\Omega }|D(\mathbf{z}_{\nu })|^{2}\,d\mathbf{x}\leqslant f(t)(||%
\mathbf{z}_{\nu }\Vert _{C(\overline{\Omega })}^{2}+1)\leqslant Cf(t)(||%
\mathbf{z}_{\nu }\Vert _{L_{2}(\Omega )}^{2}+\Vert \omega _{\nu }\Vert
_{L_{p}(\Omega )}^{2}+1)  \label{emb2}
\end{equation}%
with $f(t)\in L_{1}(0,T)$ depending only on the data $a$, $b$, $\alpha $
(independently on $\nu $) due to \eqref{eq00sec12} and \eqref{ell-1}$_{3,4}$. In the last inequality
we used
\begin{equation}
|\boldsymbol{\lambda }|^{2}\leqslant C||\mathbf{z}_{\nu }\Vert
_{L_{2}(\Omega )}^{2}\text{\quad }\;\text{ a.e. in}\;(0,T)  \label{emb3}
\end{equation}%
which holds due to \eqref{inv}, \eqref{VVV} and $\mathbf{u}_{{\nu }}\perp \mathbf{u}_{{k}}$, 
$\forall k=1,...,n.$ Applying Gronwall's Lemma to \eqref{emb2}, we
get \eqref{L2estimate}.$\hfill \;\blacksquare $

\medskip 
In the following Lemma we estimate the velocity $\mathbf{v}_{\nu }$ \ and
the vorticity $\omega _{\nu }$, independently of $\nu $.
\begin{lemma}
\label{teo4sec2} Assume that the hypothesis (\ref{eq00sec12}) hold, then we have
the estimates
\begin{align}
||\omega _{\nu }||_{L_{\infty }(0,T;L_{p}(\Omega ))}& \leqslant C, \quad
\quad \nu ||\nabla \omega _{\nu }||_{L_{2}(\Omega _{T})}^2 \leqslant C,
\label{w2} \\
||\mathbf{v}_{\nu }||_{L_{\infty }(0,T;\;W_{p}^{1}(\Omega ))}& \leqslant
C,\quad
\quad \nu ||D( \mathbf{v}_{\nu }-\mathbf{a} )||_{L_{2}(\Omega _{T})}^2 \leqslant C,
\quad \quad ||\partial _{t}\mathbf{v}_{\nu }||_{L_{2}(\Omega_{T})}\leqslant C,   \label{h2}
\end{align}%
where the constants $C$ are independent of $\nu .$
\end{lemma}

\textbf{Proof}. Let us define $B=B(\mathbf{x},t)$ as the solution of the
system 
\begin{equation*}
-\Delta B=0\quad \text{in}\;\;\Omega ,\quad \quad B\big|_{\Gamma }=b\quad 
\text{ for }t\in \left( 0,T\right).
\end{equation*}%
By Theorems 1.4 and 1.7 of Amrouche and  Rodriguez-Belido \cite{amr} the function $B$ satisfies the inequalities
\begin{equation*}
||B(\cdot ,t)||_{W_{p}^{1}(\Omega )}\leqslant C||b(\cdot ,t)||_{W_{p}^{1-%
\frac{1}{p}}(\Gamma )},\;\;\;||\partial _{t}B(\cdot ,t)||_{L_{p}(\Omega
)}\leqslant C||\partial _{t}b(\cdot ,t)||_{W_{p}^{-\frac{1}{p}}(\Gamma )}.
\end{equation*}%
According to (\ref{eq00sec12}) we
get 
\begin{equation}
\left\{ ||\partial _{t}B(\cdot ,t)||^2_{L_{p}(\Omega )}+||\bigtriangledown
B(\cdot ,t)||_{L_{p}(\Omega )}^{2}\right\} \in L_{1}(0,T).  \label{eq5sec3}
\end{equation}

Now we are able to deduce the estimates (\ref{w2}). Let us consider the
function $z_{\nu }:=\omega _{{\nu }}-B,$ which verifies the system%
\begin{equation}
\left\{ 
\begin{array}{l}
\partial _{t}z_{{\nu }}+\mathrm{div}\,(\mathbf{v_{{\nu }}\ }z_{{\nu }})={\nu 
}\Delta z_{{\nu }}+F_{{\nu }}\quad \quad \text{in }\Omega _{T}, \\ 
z_{{\nu }}|_{\Gamma _{T}}=0,\quad \quad \quad \quad z_{{\nu }}|_{t=0}=\omega
_{{0}}-B|_{t=0}\quad \quad \text{in }\Omega 
\end{array}%
\right.   \label{z}
\end{equation}%
with $F_{{\nu }}:=-\partial _{t}B-%
\mathbf{v}_{{\nu }}\cdot \nabla B.$ 
Multiplying the equality in (\ref{z}) by $q|z_{{\nu }}|_{\delta
}^{q-1}sgn_{\delta }\left( z_{{\nu }}\right) $ with $|z_{{\nu }}|_{\delta }:=%
\sqrt{z_{{\nu }}^{2}+\delta ^{2}}\,$ and $sgn_{\delta }\left( z_{{\nu }%
}\right) :=\frac{z_{{\nu }}}{\sqrt{z_{{\nu }}^{2}+\delta ^{2}}},$ we can
write%
\begin{equation}
\partial _{t}(|z_{{\nu }}|_{\delta }^{q})+\mathrm{div}(\mathbf{v}_{{\nu }%
}\,|z_{{\nu }}|_{\delta }^{q})=\big({\nu }\,\bigtriangleup z_{{\nu }}+F_{{%
\nu }}\big)\,q|z_{{\nu }}|_{\delta }^{q-1}sgn_{\delta }\left( z_{{\nu }%
}\right) .  \label{L_p}
\end{equation}%
Moreover integrating over $\Omega $ \ for $q:=p$ and taking the limit as $%
\delta \rightarrow 0^{+},$ we obtain%
\begin{align}
\frac{d}{dt}||z_{{\nu }}||_{L_{p}(\Omega )}^{p}& +{\nu }p(p-1)\int_{\Omega
}|z_{{\nu }}|^{p-2}|\nabla z_{{\nu }}|^{2}\,d\mathbf{x}  \notag \\
&\leqslant C \bigl[ ||\partial _{t}B||_{L_{p}(\Omega )} +||\bigtriangledown
B||_{L_{p}(\Omega )}||\mathbf{v}_{{\nu }}||_{C(\overline{\Omega })} \bigr]
||z_{{\nu }}||_{L_{p}(\Omega )}^{p-1} .  \label{r}
\end{align}%
From  \eqref{ell-1}$_3$, \eqref{emb1}, \eqref{emb3} we have%
\begin{equation}
||\mathbf{v}_{{\nu }}||_{C(\overline{\Omega })}\leqslant C\Vert \mathbf{v}_{{%
\nu }}\Vert _{W_{p}^{1}(\Omega )}\leqslant C(\Vert \omega _{{\nu }}\Vert
_{L_{p}(\Omega )}+||\mathbf{v}_{\nu }-\mathbf{a}\Vert _{L_{2}(\Omega
)}+||a||_{W_{p}^{1-\frac{1}{p}}(\Gamma )}),  \label{r1}
\end{equation}%
hence applying in \eqref{r}  the Young inequality, (\ref{L2estimate}), (\ref%
{eq5sec3}), (\ref{r1}) and the Bihari (Gronwall's type) inequality, we
derive the first estimate of \ (\ref{w2}).  The estimate (\ref{r}) is valid in the case when $p=2$, then
knowing  (\ref{w2})$_1$, we deduce the second inequality
of \ (\ref{w2}).   

Using (\ref{L2estimate}),  \ (\ref{w2})$_1$ and (\ref{r1}),
we obtain the first two estimates of \ (\ref{h2}). The functions $\overline{h}_{\nu }:=\partial _{t}h_{\nu }$ and $%
\overline{h}_{a}:=\partial _{t}h_{a}$ satisfy the following system
\begin{equation}
-\Delta \overline{h}_{\nu }=\,G_{\nu }  \;\;\;\text{in}\;\;\Omega \quad\text{and} 
\quad \overline{h}_{\nu }=0  \;\;\;\text{on}\;\;\Gamma 
  \qquad \text{ a.e. on  } \quad (0,T) \label{dert}
\end{equation}%
with $G_{\nu }:=\mbox{div}(\mathbf{g}_{{\nu }})$ and $\mathbf{g}_{{\nu }}=-%
\mathbf{v}_{{\nu }}\,\omega _{{\nu }}+\nu \,\nabla \omega _{{\nu }},$ such
that $\ ||\mathbf{g}_{{\nu }}||_{L_{2}(\Omega _{T})}\leqslant C$ due to  
(\ref{w2}) and (\ref{h2})$_{1}$. By Proposition 1.1, p. 12 in Girault 
and Raviart \cite{gir} and (\ref{eq00sec12}), we obtain 
\begin{equation}
||\partial _{t}\mathbf{u}_{\nu }||_{L_{2}(\Omega _{T})}\leqslant C||\mathbf{g%
}_{{\nu }}||_{L_{2}(\Omega _{T})}.  \label{dert1}
\end{equation}%
Moreover multiplying the first equation in \eqref{ICNS} by $\mathbf{u}_{i}$
and integrating over $\Omega $, we obtain%
\begin{equation}
\mathbf{A}\frac{d\boldsymbol{\lambda }}{dt}=\overline{\mathbf{f}}\quad \quad 
\text{with}\mathbf{\mathbf{\quad }}\overline{\mathbf{f}}=(\overline{f}%
_{1},...,\overline{f}_{n})^{T}\mathbf{\quad }\text{and}  \notag
\end{equation}%
\begin{equation*}
\overline{f}_{i}(t)=-\int_{\Omega }\left\{ \partial _{t}(\mathbf{u}_{\nu }+%
\mathbf{a})\cdot \mathbf{u}_{i}-\omega _{\nu }(\mathbf{v}_{\nu }\times 
\mathbf{u}_{i})\right\} \,d\mathbf{x}+\nu \int_{\Gamma }(g-2k~\mathbf{v}%
_{\nu }\cdot {\mathbf{\mathsf{s}}})(\mathbf{u}_{i}\cdot {\mathbf{\mathsf{s}}}%
)\,d\mathbf{x}-2\nu \int_{\Omega }D(\mathbf{v}_{\nu }):D(\mathbf{u}_{i})\,d%
\mathbf{x}.
\end{equation*}%
Here $\mathbf{z}\times \mathbf{u}:=z_{1}u_{2}-z_{2}u_{1}.$ Applying %
\eqref{inv},  (\ref{ell-1})$_{4}$, (\ref{w2})$_{1}$, (\ref{h2})$_{1, 2}$ and (\ref{dert1}),
we get 
$
\int_{0}^{T}|\frac{d\boldsymbol{\lambda }}{dt}|^{2}dt\leqslant C,
$
that, combining with (\ref{VVV}), (\ref{ell-1})$_{4}$, (\ref{dert1}), implies the last estimate
of (\ref{h2}).  $\hfill \;\blacksquare $

\subsection{Boundary condition on the inflow region}

\label{sec453}

The third Lemma shows that the inviscid limit preserves the boundary
condition on the inflow region $\Gamma _{T}^{-}.$ Let us define the distance
function $d(\mathbf{x}):=\min_{\mathbf{y}\in \Gamma }|\mathbf{x}-\mathbf{y}|$
from any $\mathbf{x}\in \Omega $ to $\Gamma $. Since $\Gamma \in C^{2},$ the
function $d=d(\mathbf{x})\in C^{2}(U_{\sigma _{0}}(\Gamma ))$ for a small $%
\sigma _{0}>0$ and $\bigtriangledown d=-\mathsf{n}$ on $\Gamma $. Here $%
U_{\sigma }(\Gamma ):=\left\{ \mathbf{x}\in \Omega :\ d(\mathbf{x)}<\sigma
\right\} $ is a neighborhood of $\Gamma .$

The next lemma is fundamental to show that the limit function of the
sequence $\left\{\omega _{{\nu }}\right\} $ verifies the boundary condition
on $\Gamma _{T}^{-}$ in the sense of the equality (\ref{eq7}). Here we follow the methods in Chemetov and Antontsev \cite{C}, 
Chemetov et al. \cite{CCG}.

\begin{lemma}
\label{lem6sec4} For any \underline{non-negative} test function $\psi ,$ we
have%
\begin{align}
\lim_{{\sigma }\rightarrow 0^{+}}\,\,\left( \,\overline{\lim_{{\nu }%
\rightarrow 0^{+}}}\,\,\frac{1}{{\sigma }}\int_{\Omega _{T}\cap \lbrack
\sigma <d<{2\sigma }]}|\omega _{{\nu }}-B|\;\left( \mathbf{v}_{{\nu }%
}\bigtriangledown d\right) \,\,\psi \;d\mathbf{x}dt\,\right) =0.
\label{eq35sec3}
\end{align}
\end{lemma}

\textbf{Proof}. \ Let us multiply \eqref{L_p} with $q:=1$ by a non-negative
function $\eta \in W_{2}^{1,1}(\Omega _{T})$ satisfying (\ref{psi1}), and
integrate over $\Omega .$ We get 
\begin{align*}
-\int_{\Omega _{T}}|z_{{\nu }}|_{\delta }\{\eta _{t}+\mathbf{v}_{{\nu }%
}\cdot \nabla \eta \}\;d\mathbf{x}dt& +\delta \int_{\Gamma _{T}}\,\,a\,\eta
\;d\mathbf{x}dt-\int_{\Omega }|z_{{\nu }}(0)|_{\delta }\eta (0)d\mathbf{x} \\
& \leqslant \int_{\Omega _{T}}\left\{ {\nu }|\nabla z_{{\nu }%
}| |\nabla \eta |+|F_{{\nu }}|\,\,\eta \right\} \;d%
\mathbf{x}dt.
\end{align*}%
Taking $\delta \rightarrow 0^{+}$ and using the estimates (\ref{w2})-(\ref%
{h2}), we obtain%
\begin{eqnarray}
-\int_{\Omega _{T}}\,\left( \mathbf{v}_{{\nu }}\cdot \nabla \eta \right) \
|z_{{\nu }}|\;d\mathbf{x}dt &\leqslant &\int_{\Omega }|z_{{\nu }}(0)|
\eta (0) d\mathbf{x}+\int_{\Omega _{T}}\left\{ \eta _{t}|z_{{%
\nu }}|+{\nu } |\nabla z_{{\nu }}| |\nabla \eta
|+|F_{{\nu }}|\,\,\eta \right\} \;d\mathbf{x}dt  \notag \\
&\leqslant &C\left( \int_{\Omega _{T}}(|\eta _{t}|^{2}+\sqrt{\nu }|\nabla
\eta |^{2}+\eta ^{2})\ d\mathbf{x}dt\right) ^{1/2}.  \label{eq36sec3}
\end{eqnarray}%
Let us consider the function 
\begin{equation}
{{1}_{\sigma }}(s):=\left\{ 
\begin{array}{l}
0\quad \quad \mbox{ if }s\leqslant \sigma \quad \text{and}\quad 1,\quad
\quad \mbox{ if }2\sigma \leqslant s, \\ 
\frac{s-\sigma }{\sigma },\quad \mbox{ if }\sigma <s<2\sigma ,%
\end{array}%
\right.   \label{1d}
\end{equation}%
and take $\eta :=(1-1_{\sigma }(d(\mathbf{x})))\,\psi $ in (\ref{eq36sec3}),
where $\psi $ is a non-negative test function. With the use of the
inequalities (\ref{w2})$_1$, (\ref{h2})$_1$ we deduce 
\begin{equation}
\int_{\Omega _{T}\cap \lbrack \sigma <d<{2\sigma }]}\frac{\mathbf{v}_{{\nu }%
}\cdot \nabla d}{{\sigma }}\;\psi \,|z_{{\nu }}|\,d\mathbf{x}dt\leqslant
C\int_{\Omega _{T}}\sqrt{\nu }|\nabla \eta |^{2}+(1-1{_{\sigma }}%
)^{2}(|\nabla _{\mathbf{x},t}\psi |^{2}+\psi ^{2})\;d\mathbf{x}dt.
\label{eq1122sec2}
\end{equation}%
By \eqref{emb}, the estimates (\ref{h2})$_{1,3}$ and
Corollary 9 of Simon \cite{sim}, 
\begin{equation}
\text{the set}\quad \left\{ \mathbf{v}_{{\nu }}\right\} _{\nu >0}\quad \text{%
is relatively compact in}\quad C(0,T;C^{\alpha }(\overline{\Omega })), \quad \alpha <1-2/p,
\label{compact}
\end{equation}%
having $\mathbf{v}_{\nu }\cdot \nabla d=-a$ on $\Gamma _{T}.$ Therefore
there exists $\sigma _{1}<\sigma _{0},$ independent of $\nu $, such that 
\begin{equation*}
\,\mathbf{v}_{{\nu }}\cdot \bigtriangledown d>0\qquad \text{in}\quad \cup
_{0<t<T}U_{\sigma _{1}}(\Gamma ^{-}(t)),
\end{equation*}%
where $U_{\sigma }(\Gamma ^{-}(t)):=\left\{ \mathbf{x}\in U_{\sigma }(\Gamma
):\;p(\mathbf{x})\in \Gamma ^{-}(t)\right\} $ is a neighborhood of $\Gamma
^{-}(t)$ and the projection point $p(\mathbf{x})\in \Gamma $ of $\mathbf{x}$
satisfies $|p(\mathbf{x})-\mathbf{x}|=d(\mathbf{x}).$

Due to (\ref{psi1}) and (\ref{eq1122sec2}), there exists some $\sigma
_{2}<\sigma _{1}$, such that for $2\sigma <\ \sigma _{2}$ 
\begin{align*}
0\leqslant \overline{\mathop{\lim}\limits_{{\nu }\rightarrow 0^{+}}}\,\,%
\frac{1}{\sigma }\int_{\Omega _{T}\cap \lbrack \sigma <d<{2\sigma }]}\left( 
\mathbf{v}_{{\nu }}\cdot \nabla d\right) \;\psi \,|z_{{\nu }}|\,d\mathbf{x}%
dt \leqslant C\int_{\Omega _{T}}(1-1{_{\sigma }})^{2}(|\nabla _{\mathbf{x}%
,t}\psi |^{2}+\psi ^{2})\;d\mathbf{x}dt.
\end{align*}%
Since $1_{\sigma }(d)\underset{{\sigma }\rightarrow 0^{+}}{\longrightarrow }%
1 $ in $\Omega ,$ we derive (\ref{eq35sec3})$.\hfill \;\blacksquare $

\medskip  According to (\ref{w2})-(\ref{h2}) and (\ref{compact}), there
exists a subsequence of $\nu \rightarrow 0^{+},$ satisfying (\ref{transition}%
). The pair $\{\omega ,\mathbf{v}\}$ satisfies the system (\ref{eq8})-(\ref%
{eq10}).

Now we prove that (\ref{eq7}) holds. Multiplying the parabolic equation of (\ref{nsw})
 \ by $\eta _{\sigma }:=1_{\sigma }(d(\mathbf{x}))\,\psi ,$ where $\psi $
is an arbitrary test function and integrating over $\Omega _{T},$ we obtain 
\begin{align*}
0=\{\int_{\Omega _{T}}\omega _{{\nu }}(\psi _{t}+\mathbf{v}_{{\nu }}\cdot
\nabla \psi )1{_{\sigma }}& -{\nu \ (}\nabla \omega _{{\nu }}\cdot \nabla
\eta _{\sigma })\;d\mathbf{x}dt+\int_{\Omega }\omega _{0}(\mathbf{x})\,\eta
_{\sigma }(\mathbf{x},0)\;d\mathbf{x}\,\} \\
& +\frac{1}{{\sigma }}\int_{\Omega _{T}\cap \lbrack \sigma <d<{2\sigma }%
]}\omega _{{\nu }}\,\left( \mathbf{v}_{{\nu }}\cdot \nabla d\right) \,\psi
\;d\mathbf{x}dt=J^{{\nu },\sigma }+I^{{\nu },\sigma }.
\end{align*}%
Using (\ref{transition}) and $1_{\sigma }(d)\mathop{\longrightarrow}\limits_{%
{\sigma }\rightarrow 0^{+}}1$ in $\Omega ,$ we have%
\begin{equation*}
\lim_{{\sigma }\rightarrow 0^{+}}\,\,\left( \,\lim_{{\nu }\rightarrow
0^{+}}J^{{\nu },\sigma }\right) =\int_{\Omega _{T}}\omega (\psi _{t}+\mathbf{%
v}\cdot \nabla \psi )\;d\mathbf{x}dt+\int_{\Omega }\omega _{0}\psi (\mathbf{x%
},0)\;d\mathbf{x.}
\end{equation*}%
Moreover%
\begin{equation*}
I^{\nu ,\sigma }=\{\frac{1}{\sigma }\int_{\Omega _{T}\cap \lbrack \sigma
<d<2\sigma ]}z_{{\nu }}\,\left( \mathbf{v}_{{\nu }}\cdot \nabla d\right)
\,\psi \;d\mathbf{x}dt\}+\{\frac{1}{\sigma }\int_{\Omega _{T}\cap \lbrack
\sigma <d<2\sigma ]}B\,\left( \mathbf{v}_{{\nu }}\cdot \nabla d\right)
\,\psi \;d\mathbf{x}dt\}=I_{1}^{\nu ,\sigma }+I_{2}^{\nu ,\sigma }.
\end{equation*}%
From Lemma \ref{lem6sec4}, we have $\displaystyle{\lim_{{\sigma }\rightarrow
0^{+}}\,(\,\overline{\displaystyle{\lim_{\nu \rightarrow 0^{+}}}}\,|I_{1}^{{%
\nu },\sigma }|)}=0.$

Since the function $B$ has the trace value $b$ on the boundary $\Gamma
_{T}^{-},$ the set $\left\{ \mathbf{v}_{{\nu }}\right\} $ \ is uniformly
continuous on $\overline{\Omega }_{T}$ independently on $\nu $ and $\mathbf{v%
}_{{\nu }}\cdot \bigtriangledown d=-a$ on $\Gamma _{T}$ (see the proof of
Lemma \ref{lem6sec4}), we can conclude 
\begin{equation*}
\lim_{{\sigma }\rightarrow 0^{+}}\,\,(\lim_{{\nu }\rightarrow
0^{+}}\,\,I_{2}^{{\ \nu },\sigma })=-\int_{\Gamma _{T}^{-}}\,ab\,\psi \;d%
\mathbf{x}dt.
\end{equation*}%
Therefore the pair $\{\omega ,\,\mathbf{v} \}$ satisfies the equation 
(\ref{eq7}).   Finally, 
 we can write  \eqref{NSVV} 
 in the distributional sense (for any divergence free test function 
with the zero normal component) and take the vanishing viscous limit, 
using (\ref{transition}). Theorem 2.3 in Girault, Raviart  \cite{gir}   
implies the existence of a single valued function  $p \in W^1_2 (\Omega_T )$, 
being unique  up an additive constant,
 such that the pair $\{\mathbf{v} ,\,  p \}$ fulfills  the Euler equations
\eqref{ICEE} a. e. in $\Omega_T $. Let us point out that the existence of 
$p \in W^1_2 (\Omega_T )$
can be shown by another way: knowing that
 $\mathbf{v}_{t}+  (\mathbf{v}\cdot \nabla )\mathbf{v}  \in L_2 (\Omega_T )$ 
by the regularity (\ref{transition}), and applying Theorems 3.1, 3.2 of \cite{gir}.   $\hfill \;\blacksquare $

\end{document}